\documentclass[twocolumn,10pt]{IEEEtran}
\usepackage[bottom=0.68in,top=0.6in,left=0.625in,right=0.625in]{geometry}

\DeclareUnicodeCharacter{FF09}{ }
\DeclareUnicodeCharacter{FF08}{  }
\DeclareUnicodeCharacter{FF0C}{    }
\DeclareUnicodeCharacter{FF1B}{     }
\DeclareUnicodeCharacter{FF1A}{      }

\ifCLASSINFOpdf
  \usepackage[pdftex]{graphicx}
  \graphicspath{{../pdf/}{../jpeg/}}
  \DeclareGraphicsExtensions{.pdf,.jpeg,.png}
\else
  \usepackage[dvips]{graphicx}
  \graphicspath{{../eps/}}
  \DeclareGraphicsExtensions{.eps}
\fi

\usepackage{epstopdf}
\usepackage{multirow}
\usepackage{float}
\usepackage[cmex10]{amsmath}
\usepackage {amssymb}
\usepackage{graphicx}
\usepackage{array}
\usepackage{mdwmath}
\usepackage{mdwtab}
\usepackage{eqparbox}
\usepackage{commath}
\usepackage{nomencl}
\usepackage{cite}
\usepackage{algorithm}
\usepackage{algorithmic}
\usepackage{mathtools}
\usepackage{makeidx}
\usepackage{ifthen}
\usepackage{subfigure}
\usepackage{caption}
\usepackage{color}
\usepackage{multirow,array}
\usepackage{pdflscape}
\usepackage{bm,upgreek}

\makenomenclature

\renewcommand{\nomgroup}[1]{%
\ifthenelse{\equal{#1}{I}}{\item[\textbf{Indices}]}{%
\ifthenelse{\equal{#1}{A}}{\item[\textbf{Abbreviations}]}{%
\ifthenelse{\equal{#1}{V}}{\item[\textbf{Variables}]}{%
\ifthenelse{\equal{#1}{P}}{\item[\textbf{Parameters and Constants}]}{%
}
}
}
}
}

\ifCLASSINFOpdf
  \usepackage[pdftex]{graphicx}
  \graphicspath{{../pdf/}{../jpeg/}}
  \DeclareGraphicsExtensions{.pdf,.jpeg,.png}
\else
  \usepackage[dvips]{graphicx}
  \graphicspath{{../eps/}}
  \DeclareGraphicsExtensions{}
\fi
\usepackage{epstopdf}

\usepackage{multirow,tikz,enumitem}
\usepackage{xspace,amsmath,amsfonts,amssymb,tikz,multirow, pgfplots,pbox}
\usepackage{caption}
\usepackage{xspace}
\usepackage{latexsym}
\usepackage{xcolor}
\usepackage{graphicx}
\usepackage{enumitem}
\usepackage{comment,verbatim}
\usepackage{array}
\usepackage{booktabs}
\usepackage{cite}
\usepackage{fancyhdr}

\usepackage{comment}
\usepackage{epstopdf}
\usepackage{float}

\newcommand{\beq}{\begin{equation}}
\newcommand{\eeq}{\end{equation}}
\newcommand{\beqn}{\begin{eqnarray}}
\newcommand{\eeqn}{\end{eqnarray}}
\newcommand{\beqno}{\begin{eqnarray*}}
\newcommand{\eeqno}{\end{eqnarray*}}
\newcommand{\bma}{\begin{displaymath}}
\newcommand{\ema}{\end{displaymath}}
\newcommand{\bnu}{\begin{enumerate}}
\newcommand{\enu}{\end{enumerate}}
\newcommand{\bce}{\begin{center}}
\newcommand{\ece}{\end{center}}
\newcommand{\btb}{\begin{tabular}}
\newcommand{\etb}{\end{tabular}}

\def\bx{{\mathbf{x}}}

\def\bu{{\mathbf{u}}}
\def\by{{\mathbf{y}}}

\allowdisplaybreaks[4]

\begin{document}

\title{Two-Stage Distributionally Robust Edge Node Placement  Under Endogenous  Demand Uncertainty}
\author{ \IEEEauthorblockN{Jiaming Cheng\IEEEauthorrefmark{2},Duong~Thuy~Anh~Nguyen\IEEEauthorrefmark{3},  Duong Tung Nguyen\IEEEauthorrefmark{2}\\
\IEEEauthorblockA{\IEEEauthorrefmark{2}Arizona State University, Tempe, AZ 85281, USA,
\{dtnguy52, duongnt\}@asu.edu} \\
\IEEEauthorblockA{\IEEEauthorrefmark{3}University of British Columbia, Vancouver, BC V6T1Z4, Canada, jiaming@ece.ubc.ca} \vspace{-0.57em}}}
\maketitle

\begin{abstract}
Edge computing (EC) promises to deliver low-latency and ubiquitous computation to numerous devices at the network edge. This paper aims to jointly optimize edge node (EN) placement and resource allocation for an EC platform, considering demand uncertainty. Diverging from existing approaches treating uncertainties as exogenous, we propose a novel two-stage decision-dependent distributionally robust optimization (DRO) framework to effectively capture the interdependence between  EN placement decisions and uncertain demands. The first stage involves making EN placement decisions, while the second stage optimizes resource allocation after uncertainty revelation. We present an exact mixed-integer linear program reformulation for solving the underlying ``min-max-min" two-stage model. We further introduce a valid inequality method to enhance computational efficiency, especially for large-scale networks. Extensive numerical experiments demonstrate the benefits of considering endogenous uncertainties and the advantages of the proposed model and approach.
\end{abstract}


\printnomenclature

\section{Introduction}
\label{intro}
The widespread adoption of mobile devices and applications has caused an unprecedented surge in mobile data traffic. Moreover, with the advent of new services such as augmented/virtual reality, manufacturing automation, and autonomous driving, the demand for innovative solutions that can meet their challenging requirements has become imperative. In response, edge computing (EC) has emerged as a vital computing paradigm, complementing traditional cloud computing to provide enhanced user experiences and support a wide array of low-latency and highly reliable Internet of Things (IoT) applications \cite{wshi16,Antonino2019}. However, as the utilization of edge resources continues to grow, it poses significant challenges to existing network operations. The intermittent nature and stringent service requirements of EC, combined with system uncertainties, create substantial obstacles to network management and optimization. 

The performance and reliability of EC systems are susceptible to various uncertainties arising  from multiple sources, such as extreme weather conditions, fluctuating resource demands, traffic spikes, user mobility, and changes in application performance and user behavior. Moreover, the increasing complexity and diversity of man-made attacks and cyber threats, including insider attacks, cyberattacks, and malware attacks,  introduce additional uncertainties and risks to EC systems \cite{Moglen23}. Indeed, various aspects of edge network operations under uncertainties have been studied, including computational resource allocation under demand uncertainty \cite{Duong_iot,ouyang19,netsoft23,wang22}, resilience network designs against EN failures \cite{ella2022,resilientEC,fhe22}, network softwarization against risks \cite{li23,dche20,pzha18,symeon20}, economic analysis under price uncertainty \cite{Lei11,niyato15,Dejan19,Ivashko20}, and  market interaction under time-varying communication network \cite{wiopt23,DO_SO,ella_TCNS}.

Incorporating demand uncertainty is of utmost importance in the long-term investment and operation plan of an EC platform when selecting potential locations for edge resource installation. The platform must make edge node (EN) placement decisions based on incomplete information about future demand. Ignoring uncertain demand may lead to frequent over-provisioning or under-provisioning of resources. Over-provisioning can result in the wastage of resources and unnecessary high provisioning costs while under-provisioning may lead to degraded service quality and unmet demand.

To address this challenge, we propose a distributionally robust model designed to address the EN placement and resource allocation problem for a budget-constrained EC platform. The strategic placement of ENs plays a crucial role in determining their proximity to users and their ability to serve specific areas. Thus, user demand patterns, data traffic, and service requests are significantly influenced by the chosen EN deployment locations. As a result, the EN placement decision directly impacts the actual demand realization. Our proposed model aims to optimize the EN placement decision under decision-dependent demand uncertainty, with the goal of minimizing costs while enhancing the quality of service (QoS) in terms of latency and unmet demand. Effectively managing uncertainties is a key enabler in achieving consistent performance, reliability, and a superior user experience in EC. 

Many efforts within the realm of optimization under uncertainty have been developed for EC, with stochastic optimization (SO) and robust optimization (RO) being the two main approaches. SO typically assumes complete knowledge of the underlying uncertainty distribution and requires access to a large number of samples drawn from this true distribution. However,  this assumption may be demanding in practice, and limited information can lead to misspecification of the distribution \cite{SO}. On the other hand, RO adopts uncertainty sets, deterministic representations of uncertain parameters, simplifying the model and improving computational tractability  \cite{RObook,RO_15}. Nevertheless, RO can be overly conservative at times, potentially leading to suboptimal system performance. The distributionally robust optimization (DRO) approach strikes a balance between SO and RO \cite{kuhn_dro,delage10}. It optimizes decisions with respect to worst-case distribution within a predefined ambiguity set, achieving a favorable trade-off between optimality and robustness. 

While DRO has been effectively applied in planning and operation problems in cloud/EC \cite{li23,Li20,DRO_zhang,Zhil23}, one aspect often neglected in the literature is the interdependence between decisions and uncertainties. Specifically, the placement of ENs in an area's neighborhood has a positive impact on demand. It boosts user confidence, especially when they request services with stringent delay requirements. Thus, the platform may expect an increased mean of demand. Moreover, with the increasing number of ENs in an area's neighborhood, higher user confidence also leads to decreased demand variance. As users become more confident in the reliability and availability of edge resources, their demand patterns tend to become more consistent and predictable. This reduced variability indicates that users exhibit a more stable and reliable demand, enabling better resource planning and management. Thus, it is crucial for the platform to consider the influence of its decisions on future demands. Consequently, the platform can proactively optimize its decisions to control the uncertainty set. Surprisingly, this critical problem has been largely overlooked in the existing literature. It presents a fundamental and unresolved challenge in optimizing EN placement decisions. Unfortunately, the existing research in computer networking lacks the necessary tools and techniques to tackle this problem effectively. 

\noindent \textbf{Contribution:} This paper seeks to bridge the gap by proposing a novel two-stage DRO framework with a decision-dependent moment-based ambiguity set for optimal EN placement. Unlike conventional DRO approaches that use exogenous ambiguity sets, our proposed model incorporates an endogenous ambiguity set, which captures the interdependence between the first and second moments of demand and the placement decisions. To the best of our knowledge, we are the first to consider this two-stage distributionally robust EN placement model that explicitly accounts for this decision-dependent demand uncertainty. However, incorporating the interdependence between uncertainties and decisions increases complexity, resulting in a large-scale non-linear optimization problem with numerous bilinear and trilinear terms. To tackle the challenging problem, we first develop an efficient and exact reformulation, termed \textbf{Exact OPT-Placement}. This reformulation is achieved through \textit{3-step} transformations that convert the problem into a Mixed Integer Linear Programming (MILP) form, which can be solved efficiently using widely available solvers (e.g., Gurobi and Mosek). We further introduce an improved algorithm that generates feasibility cuts to strengthen the proposed algorithm and speed up the computation. To substantiate the effectiveness of our approach, extensive simulations have been conducted, demonstrating the efficiency of the proposed scheme in comparison to several baseline models. Additionally, we have performed sensitivity analyses to evaluate the impact of crucial system parameters on the overall system performance. 

\allowdisplaybreaks[4]

\section{System Model and Problem Formulation}
\label{system}
In this section, we present the DRO model for the EN placement and resource allocation problem for a budget-constrained EC platform. The main objective of the platform is to optimize the EN placement decision under endogenous demand uncertainty, aiming to minimize costs while enhancing QoS.

\subsection{System Model}
We consider an EC platform that manages a set $\mathcal{J}$ of $J$ potential candidate locations for EN installation and provides edge resources to users in a set $\mathcal{I}$ of $I$ areas, each represented by an access point (AP). The AP and EN indices are denoted by $i$ and $j$, respectively. The size of ENs can vary significantly, and each EN may comprise one or multiple edge servers. For simplicity, we consider only computing resources, and the resource capacity at EN to be placed at location $j$ is denoted by $C_j$. It is straightforward to extend our model to consider the sizing decision for each EN. Given the diverse range of IoT services with varying requirements, edge servers are responsible for hosting different types of IoT applications to serve these workloads effectively. The platform optimizes the long-term EN placement in the initial stage, maintaining this configuration unchanged for an extended period. 

The placement decision for an EN at location $j\!\in\! \mathcal{J}$ is denoted by a binary variable $y_j \!\in\! \{0,1\}$. Specifically, $y_j$ takes the value $1$ if an EN is installed at location $j$ and $0$ otherwise. Additionally, this placement decision incurs an EN placement cost of $f_j$. The objective of the platform is to identify the optimal set of locations for efficient EN placement while adhering to the budgetary constraints imposed by the investment budget $B$. The network delay between AP $i$ and EN $j$ is $d_{i,j}$. To minimize network delay, the demand in each area should ideally be served by its closest EN. However, each EN has a limited capacity. Therefore, given the first-stage decision, the platform needs to optimally allocate the actual demand to different ENs, considering the edge resource capacity constraints and the diverse geographical locations of the ENs to ensure high QoS while reducing costs. Since the exact demand is unknown to the platform in the first stage, a portion of the workload may be dropped. Let $x_{i,j}$ denote workload from area $i$ allocated to EN $j$ and $u_i$ be the unmet demand from area $i$. User requests from each area $i$ must be either served by some ENs or dropped (i.e., counted as unmet demand $u_i$), and the penalty for each unit of unmet demand is denoted by $s_i$. 

The resource demand in area $i$ is denoted by $\lambda_i$. The demand $\lambda_i$ in each area $i \in \mathcal{I}$ exhibits inherent uncertainties and can vary based on the EN placement decisions. Indeed, the decisions regarding EN placement in adjacent to user-populated areas can significantly influence users' demand. The presence of more ENs, along with increased resource availability and reduced network delay, contributes to higher user confidence. Consequently, a larger pool of potential customers is attracted to utilize the available edge resources, resulting in higher demand. The platform's objective is to minimize the variation of demand (forecast error) that deviates from the initially projected demand, as it plays a crucial role in delivering an exceptional user experience. Therefore, accurately capturing and modeling this uncertain demand, which not only varies over time but also changes based on the first-stage EN placement decision, is of paramount importance. Developing efficient models for optimal placement of ENs that can effectively mitigate the uncertainty associated with demand and enhance the robustness of the system is a critical focus of this work. 

In the following, we present a deterministic model for the EN placement and workload allocation problem, followed by two uncertainty models addressing demand uncertainties: one for exogenous uncertainties and another considering the endogeneity between placement decisions and uncertainties.


\subsection{Deterministic problem formulation}
\label{deterministic}
The EC platform aims to minimize the total EN placement cost while enhancing user experience by reducing the overall network delay and unmet demand. Without uncertainty consideration, the deterministic problem for EN placement and workload allocation can be formulated as follows:
\begin{subequations}
\label{DET}
\begin{align}
     \textbf{DET:} ~~ &\min_{\by,\bx,\bu}  \:  \sum_{j \in \mathcal{J}} f_{j} y_{j}  + \rho \sum_{i \in \mathcal{I}} \sum_{j \in \mathcal{J}}  d_{i,j} x_{i,j} + \sum_{i \in \mathcal{I}} s_i u_i,  \label{eq-DETobj} \\
    & ~~ \text{s.t.} ~~ (\ref{budget}) - (\ref{var_constr1}).
\end{align}
\end{subequations} 
The first term in the objective function \eqref{eq-DETobj} represents the total EN placement cost, while the second and third term captures the delay penalty and unmet demand penalty from the workload allocation decision. $\rho$ is the delay penalty parameter controlled by the platform. A higher value of $\rho$ signifies that the platform prioritizes reducing delays over minimizing unmet demand. Constraints are summarized as follows: \\
\textbf{\textit{Budget constraint:}} 
The total expense for EN placement decisions should not exceed the restricted investment budget $B$:\vspace{-0.2cm}
\begin{align}
    \sum_{j \in \mathcal{J}} f_{j} y_{j} \leq B.  \label{budget}
\end{align}
\textbf{\textit{Reliability constraint:}} 
To enhance service reliability, the platform may opt to place the EN on a minimum of $K^{\sf min}$ locations, proactively considering the potential for unexpected EN failures. Hence, we have:\vspace{-0.2cm}
\begin{align} 
    \sum_{j \in \mathcal{J}} y_{j} \geq K^{\sf min}.  \label{reliability}
\end{align}
\textbf{\textit{Workload allocation constraints:}} 
The demand from each area must be either served by some ENs ($x_{i,j}$) or dropped ($u_i$):\vspace{-0.2cm}
\begin{align}
    u_i + \sum_{j \in \mathcal{J}} x_{i,j} = \lambda_i, ~ \forall i. \label{supply_demand}
\end{align}
\textbf{\textit{Capacity constraints:}} 
We assume that the capacity of each EN is pre-divided for each area. By allocating specific capacities to each area, the system can achieve an equitable distribution of resources among customers, preventing any particular area from being disproportionately overloaded or underserved. Let $C_{i,j}$ denote the resource capacity of each EN $j$ dedicated to area $i$.  The computing resources required to serve the workload from area $i$, assigned to  EN $j$, must not exceed the pre-allocated capacity $C_{i,j}$ designated for area $i$:\vspace{-0.2cm}
\begin{align}
    0 \leq x_{i,j} \leq C_{i,j} y_j, ~ \forall i,j. \label{resource_cap}
\end{align}
\textbf{\textit{Delay constraints:}} 
To maintain a satisfactory user experience, the platform may impose an average network delay threshold $\Delta_{i}$ for area $i$ based on the proportion of workload from area $i$ allocated to EN $j$, represented by $\frac{x_{i,j}}{\lambda_i}$. \vspace{-0.2cm}
\begin{align}
    \sum_{j \in \mathcal{J}} \frac{d_{i,j} x_{i,j}}{\lambda_i} \leq \Delta_{i}  , ~ \forall i. \label{delay}
\end{align}

\noindent \textbf{\textit{Decision variables:}} Decision variables include the EN placement decision $\by$, workload allocation $\bx$ and unmet demand $\bu$. 
\begin{align}
    \by \in \{0,1\}^{J} , ~~\bu \in \mathbb{R}_{+}^{I}, ~~ \bx \in \mathbb{R}_{+}^{I \times J} \label{var_constr1}.
\end{align}

\subsection{Uncertainty Modeling}
In the deterministic model \textbf{DET}, the demand $\lambda_i$ for each geographical area is assumed to be precisely known at the decision-making stage. This implies that the platform can predict the actual demand, which subsequently serves as an input to the \textbf{DET} model. Nevertheless, accurately predicting the exact demand for each area is often challenging at the time of decision-making. Consequently, addressing how to effectively incorporate uncertainties into the decision-making process becomes a critical and complex task. In this study, we explore an alternative approach to enhance the robustness of our model against the model uncertainty and limited information, known as the DRO approach. The core idea of DRO revolves around considering an ambiguity set that encompasses a range of possible distributions consistent with the available information. Specifically, we focus on scenarios where only the mean and variance of the demand distribution are provided.

We assume that the true distribution of demand originates from a set of possible distributions, where the random demand at each area $i \in \mathcal{I}$ can take values from a finite support set $\Xi = \{\xi_1, \xi_2, \dots, \xi_N\}$ with unknown probabilities ($p_{i,1}, p_{i,2},\dots, p_{i,N}$). Drawing upon the work in \cite{delage10}, we utilize a moment-based ambiguity set that considers information about the distribution's support and confidence region. The intuition behind this ambiguity set is to maintain the forecast error for the first and second moments of demand within specified thresholds $\Gamma_i^{\mu}$ and $\underline{\Gamma}_i^{\mu}$, $\bar{\Gamma}_i^{\sigma}$ ($0 \leq \underline{\Gamma}_i^{\mu} \leq 1 \leq \bar{\Gamma}_i^{\sigma}$), respectively, ensuring robustness against uncertainties in the prediction values. These predefined parameters play a crucial role in determining the robustness of the ambiguity set for each specific area $i$. 


\subsubsection{Exogenous Stochastic Demand} 
Let $\bar{\mu}_i$ and $\bar{\sigma}_i^2$ represent the empirical mean and variance of demand $\lambda_i$ at area $i$, respectively. When the demand is independent of the system decision, the ambiguity set for the exogenous stochastic demand can be presented as follows: \vspace{-0.2cm}
\begin{subequations}
\label{diu}
\begin{align}
    \mathcal{U}_1 (y) = \bigg\{  \{ p_{i} \}_{i \in \mathcal{I}}: ~  p_i \in \mathbb{R}_{+}^{N}, ~ \sum_{n = 1}^{N} p_{i,n} = 1, &~\forall i, \label{diu_1} \\
    \abs{\sum_{n = 1}^{N }p_{i,n} \xi_{n} - \bar{\mu}_i} \leq \Gamma_i^{\mu}, & ~\forall i, \label{diu_2}\\
    \big( \bar{\sigma}_i^2 + \bar{\mu}_i^2 \big) \underline{\Gamma}^{\sigma}_{i} \leq \sum_{n = 1}^{N } p_{i,n} \xi_{n}^2 \leq \big( \bar{\sigma}_i^2 + \bar{\mu}_i^2 \big)  \bar{\Gamma}_{i}^{\sigma}, &~\forall i \label{diu_3} ~ \bigg\}. 
\end{align}
\end{subequations}
Here, constraint (\ref{diu_1}) ensures that the probabilities across all areas within the support set sum up to $1$. In (\ref{diu_2}), the true mean of demand is constrained to lie within an $L_1$-distance $\Gamma_i^\mu$ from the empirical mean $\bar{\mu}_i$. Constraint (\ref{diu_3}) implies that the actual value of the second moment of demand must fall within the interval $[ (\bar{\sigma}_i^2 + \bar{\mu}_i^2 \big) \underline{\Gamma}_{i}^{\sigma} , (\bar{\sigma}_i^2 + \bar{\mu}_i^2 ) \bar{\Gamma}_{i}^{\sigma} ]$.

\subsubsection{Endogenous Stochastic Demand}
The presence of ENs in the neighborhood of each area $i$ has a substantial impact on increasing the demand for that area. This effect can be attributed to the close proximity of ENs, ensuring higher availability of resources and the overall QoS, especially for time-sensitive applications. As a result, user confidence is heightened, leading to a rise in the mean of demand and a decrease in demand variance since customers' demand patterns become more consistent and predictable, offering significant advantages for resource planning and management.

To capture the interdependence between placement decisions $y$ and the demand uncertainty, we assume that the demand at each area $i$ is represented by a random variable $\lambda_i(y)$, with its distribution dependent on the EN placement decisions $y$, having mean $\mu_i(y)$ and variance $\sigma_i^{2}(y)$. The endogenous decision-dependent ambiguity set $\mathcal{U}(y)$, defining bounds for the first two moments of the demand distribution as functions of $y$, \cite{luo20, Yu22}, is expressed as follows:\vspace{-0.2cm}
\begin{subequations}
\label{ddu}
\begin{align}
    & \mathcal{U}_2 (y) = \bigg\{  \{ p_{i} \}_{i \in \mathcal{I}}: ~  p_i \in \mathbb{R}_{+}^{N}, ~ \sum_{n = 1}^{N} p_{i,n} = 1, ~\forall i, \label{ddu_1} \\
    & ~~~~~~~~~~~~~~~~~~~~~~ \abs{\sum_{n = 1}^{N }p_{i,n} \xi_{n} - \mu_i (y)} \leq \Gamma_i^{\mu}, ~ \forall i, \label{ddu_2}\\
    & \big[ \sigma_i^{2}(y) \!+\! (\mu_i({y}))^2 \big] \underline{\Gamma}^{\sigma}_{i} \!\leq\! \sum_{n = 1}^{N } p_{i,n} \xi_{n}^2 \!\leq\! \big[ \sigma_i^{2}(y) \!+\! (\mu_i({y}))^2 \big] \bar{\Gamma}_{i}^{\sigma}, \!\forall i \label{ddu_3} \bigg\}\!. 
\end{align}
\end{subequations}
Similar to (\ref{diu}), (\ref{ddu}) also imposes linear restrictions on the first and second moments of the demand distribution, ensuring they remain within predefined thresholds. However, it is crucial to emphasize that in contrast to (\ref{diu}), the placement decision $y$ has a substantial impact on both the mean ($\mu_i(y)$) and variance ($\sigma^2_i(y)$) of the demand. We model the mean and variance of demand as affine functions of the decisions $y$, which can be expressed as:
\begin{subequations}
\label{moment_def}
\begin{align}
    & \mu_i(y) = \bar{\mu}_i \Bigg( 1 + \sum_{j \in \mathcal{J}} \Psi_{i,j}^{\mu} y_{j} \Bigg), \\
    & \sigma_i^2(y) = \max \Bigg\{  \bar{\sigma}_i^{2} \Bigg( 1 - \sum_{j \in \mathcal{J}} \Psi_{i,j}^{\sigma} y_{j} \Bigg), (\sigma_i^{LB})^2  \Bigg\}.
\end{align}
\end{subequations}
As previously mentioned, an increase in resource availability is generally correlated with an increase in the mean and a decrease in the variance of demand \cite{luo20, yu2022multistage}. To capture the influence of placing EN $j$ near area $i$ on the mean and variance of demand, we employ two parameters: $\Psi_{i,j}^{\mu}$ and $\Psi_{i,j}^{\sigma}$. These parameters take values from the interval $[0,1]$ and are specifically designed to reflect the varying influence of different locations. Closer locations have higher impacts on demand's first and second moments, while areas farther away exert less effect. Notably, when $\Psi_{i,j}^{\mu}$ and $\Psi_{i,j}^{\sigma}$ are set to $0$ in (\ref{moment_def}), the ambiguity set reduces to an exogenous form, as seen in (\ref{diu}). This corresponds to the traditional DRO problem, where demand uncertainty is treated independently of the placement decisions of ENs.

Furthermore, when an EN is placed near area $i$, the demand increases from the forecast demand estimate $\bar{\mu}_i$. The highest demand variance in area $i$ occurs when there are no available ENs in its neighborhood. However, due to the inherent nature of the demand uncertainty, the variance cannot be less than a lower bound value, denoted by $(\sigma_i^{LB})^2$. In situations where the platform possesses perfect knowledge of the first and second moments of demand, the parameters $\Gamma_i^{\mu} = 0$ and $\underline{\Gamma}^{\sigma}_{i} = \bar{\Gamma}_{i}^{\sigma} = 1$ can be set. In this case, the problem is reduced to \textbf{DET}, where uncertainty is eliminated, and the solution relies on known values of the forecast mean and demand variance.
\vspace{-0.2cm}
\subsection{Two-stage DRO with endogenous stochastic demand}
For any demand distribution $p=\{p_i\}_{i\in\mathcal{I}} \in \mathcal{U}(y)$, where $\mathcal{U}(y)$ is the ambiguity set defined in \eqref{ddu} with \eqref{diu} as a special case, the two-stage decision-dependent DRO problem of the EC platform for EN placement and resource allocation can be formulated as follows:
\begin{subequations}
\label{dro}
\begin{align}
    & \!\!(\mathcal{P}_1) ~ \min_{\by} \!\sum_{j}\! f_j y_j \!+\!\! \max_{p \in \mathcal{U}(y)} \!\min_{\bx, \bu} \mathbb{E}_{p}\!\bigg[\!\rho \! \sum_{i,j}\! d_{i,j} x_{i,j} \!+\!\! \sum_{i}\! s_i u_i\! \bigg] \!\!\label{dro_obj}\\
    &  \text{s.t.} ~~  (\ref{budget}),  (\ref{reliability}),  (\ref{var_constr1})  \nonumber \\
    &  \qquad \Omega_2 (y,\lambda) = \bigg\{ ~ 0 \leq x_{i,j} \leq C_{i,j} y_j, ~ \forall i,j \\
    & \qquad\qquad\qquad\qquad  u_i + \sum_{j} x_{i,j} = \lambda_i(y) , ~ \forall i \\
    & \qquad\qquad\qquad\qquad \sum_{j} d_{i,j} x_{i,j} \leq \Delta_{i} \lambda_i (y) , ~ \forall i \bigg\} .
\end{align}
\end{subequations}
The proposed two-stage DRO model addresses a trilevel min-max-min optimization problem. The first stage minimizes EN placement costs before revealing demand uncertainties. In the second stage, the model considers worst-case demand realization within an ambiguity set, aiming to minimize the expected operational penalty after demand uncertainties are realized. Importantly, the two stages are interconnected as EN placement decisions directly impact actual demand. This differs from the two-stage DRO framework with exogenous stochastic demand, where actual demand is treated as an independent instance within the ambiguity set $\mathcal{U}_1$ representing possible demand scenarios. The proposed model, however, considers actual demand after EN placement decisions have been made in the first stage.

In our problem, we assume that the unmet demand penalty ($s_i$) is larger than the network delay cost ($\rho d_{i,j}$), i.e., $s_i > \rho d_{i,j}$.  This means that the EC platform prioritizes meeting demand requirements over reducing latency in the network. 

\vspace{-0.4cm}
\section{Solution Approach}
\label{sol}
In this section, we propose two methods to solve the two-stage DRO model for EN placement and resource allocation with endogenous stochastic demand. The problem has a two-stage structure, leading to a tri-level optimization problem with a strong interdependence between uncertainty and decisions. Even a seemingly straightforward formulation employing linear programming (LP) in both stages can be proven to be NP-hard \cite{RObook}. Additionally, the uncertainty parameter $\lambda(y)$ is defined as an affine function of the decision variable $y$ within the ambiguity set, leading to substantial complexities in solving problem $(\mathcal{P}_1)$. The modification introduces significant changes to the size and structure of the set, making the reformulation of the DRO problem particularly challenging, and the presence of bilinear terms further obstructs direct solvability. 

To address these challenges, we provide an exact reformulation for problem $(\mathcal{P}_1)$, enabling us to solve the resulting MILP problem using off-the-shelf solvers. Furthermore, we propose an improved reformulation that enhances scalability by leveraging extreme rays of the feasible region to generate efficient cuts. The improved approach can accelerate computation, particularly as the network size increases.

\vspace{-0.4cm}
\subsection{Exact Monolithic Reformulation}
\label{monolithic_reform}
To reformulate the min-max-min DRO model ($\mathcal{P}_1$) into a more tractable form, we employ a three-step process to derive a single-level, monolithic representation. In \textbf{Step 1}, given EN placement decisions $y$, we derive the dual problem of the innermost minimization problem and obtain a closed-form expression for its optimal objective value. Subsequently, the inner max-min framework simplifies to a max-max formulation, effectively simplifying it to a maximization problem. In \textbf{Step 2}, we apply the duality theorem to recast the obtained min-max model as an equivalent single-level mixed-integer nonlinear program (MINLP). \textbf{Step 3} leverages the McCormick linearization techniques \cite{BigM} to attain a single-level MILP reformulation of the original tri-level problem, which can be directly solved using off-the-shelf solvers. Below are the specific details:

\textbf{\textit{Step 1:}} The structure of the dual problem allows us to decompose the inner problem based on each area $i$. Denoting the inner level problem as $g(y,\lambda)$, we express $g(y,\lambda)=\sum_{i \in \mathcal{I}}g_i(y,\lambda)$, where $g_i(y,\lambda)$ for all $i \in \mathcal{I}$ is given as follows
\begin{subequations}
\label{primal_1}
\begin{align}
    &   g_i(y,\lambda)  =  \min_{\bx, \bu} \rho \sum_{j} d_{i,j} x_{i,j} + s_i u_i\\
    &  \text{s.t.} ~~ x_{i,j} \leq C_{i,j} y_j, ~ \forall j  \qquad\qquad ~(v_{i,j}) \label{mono_1}\\
    &\qquad   u_i + \sum_{j} x_{i,j} = \lambda_i(y) ,  ~\qquad (\alpha_i) \\
    & \qquad  \sum_{j } d_{i,j} x_{i,j} \leq \Delta_{i} \lambda_i (y), ~ ~~~~ (\beta_i) \label{mono_3} 
\end{align}
\end{subequations}
where $v_{i,j}$, $\alpha_i$ and $\beta_i$ are the dual variables associated with constraints (\ref{mono_1})-(\ref{mono_3}), respectively. Consequently, the dual problem of $g_i(y,\lambda)$, for all $i \in \mathcal{I}$, is as follows:
\begin{subequations}
\label{dual_problem}
\begin{align}
    & \max_{v_{i,j},\alpha_i, \beta_i} ~~  \sum_{j} C_{i,j} y_j v_{i,j} + \big[ \alpha_i + \beta_i \Delta_{i} \big] \lambda_i(y) \label{dual_problem-obj}\\
    &  ~~\text{s.t.} ~~  v_{i,j} + \alpha_i + \beta_i d_{i,j} \leq \rho d_{i,j}, ~\forall j \label{eq-constdual}\\
    & ~~\qquad  \alpha_i \leq s_i, ~ \beta_i \leq 0; ~~  v_{i,j} \leq 0, ~ \forall j.
\end{align}
\end{subequations}
Subsequently, our goal is to derive a closed-form expression for the optimal objective value of the dual problem associated with each inner problem, considering the extreme points and rays of the feasible region. We examine the following two cases:


\noindent
\textbf{Case 1 ($\alpha_i = s_i$):} From \eqref{eq-constdual}, we have $v_{i,j} \leq \rho d_{i,j} - s_i - \beta_i d_{i,j}$. As $v_{i,j} \leq 0$, the extreme point of $v_{i,j}$ can occur at either $v_{i,j} = 0$ or $v_{i,j} = \rho d_{i,j} - s_i - \beta_i d_{i,j}$ if $\rho d_{i,j} - s_i - \beta_i d_{i,j} < 0$. We consider the following two cases for the value of $\beta_i \le 0$:

(i) If $\rho - \frac{s_i}{d_{i,j}}<\beta_i \le 0$ then $v_{i,j} \le \rho d_{i,j} - s_i<0$ due to the assumption $s_i > \rho d_{i,j}$. The constraint $v_{i,j} \le 0$ becomes redundant and $v_{i,j} = \rho d_{i,j} - s_i$ is the extreme point. The optimal value for the objective function is
\begin{align}\label{optimal_obj_case1(i)}
s_i \lambda_i(y) + \sum_{j} C_{i,j} y_j (\rho d_{i,j} - s_i), ~ \forall i.
\end{align}

(ii) If $\beta_i \leq \rho - \frac{s_i}{d_{i,j}} < 0$, the inequality $\rho d_{i,j} - s_i - \beta_i d_{i,j} > 0$ holds true. As a result, the constraint $v_{i,j} \le \rho d_{i,j} - s_i - \beta_i d_{i,j}$ becomes redundant and $v_{i,j} = 0$ represents the extreme point. By letting $d_{i}^{\sf min} = \min_{j'\in\mathcal{J}} d_{i,j'}$ for all $i$, we find that $\beta_i = \rho - \frac{s_i}{d_{i}^{\sf min}}$ is an extreme point. The objective value is
\begin{align}\label{optimal_obj_case1(ii)}
\Bigg[ s_i + \Bigg( \rho - \frac{s_i}{d_{i}^{\sf min}}\Bigg) \Delta_i \Bigg] \lambda_i(y) , ~\forall i.
\end{align}


\noindent
\textbf{Case 2 ($\alpha_i < s_i$):} For all ENs $j \in \mathcal{J}$, $v_{i,j}$ reaches its extreme point at either $v_{i,j} = 0$ or $v_{i,j} = \rho d_{i,j} - \alpha_i - \beta_i d_{i,j}$. We proceed to analyze the following two scenarios: 

(i) If $v_{i,j} = 0$ for some $j\in \mathcal{J}$, it must hold that $\rho d_{i,j} - \alpha_i - \beta_i d_{i,j}\ge 0$, i.e., $\alpha_i \le d_{i,j} (\rho - \beta_i)$. Thus, we aim to find extreme points for $\beta_i$ such that
\begin{align}
    s_i & > \bigg\{ \max_{\beta_i} (\rho - \beta_i) d_{i,j}, \forall j, ~~ \text{s.t} ~~ \beta_i \leq 0\bigg\}.\label{constr2}
\end{align}
Notably, $(\rho - \beta_i) d_{i,j}>s_i$ when $\beta_i \to -\infty, \forall i$. Thus, $\alpha_i = \rho d_{i,j}$ and $\beta_i = 0$ represent the extreme points.

(ii) If the extreme point is $v_{i,j} = \rho d_{i,j} - \alpha_i - \beta_i d_{i,j}$, it implies that the constraint $v_{i,j} \le \rho d_{i,j} - \alpha_i - \beta_i d_{i,j}$ is binding, i.e.,
\begin{align}\label{constr1}
\rho d_{i,j} - \alpha_i - \beta_i d_{i,j} \le 0.   
\end{align}
Since $\beta_i \leq 0$, $\beta_i = 0$ represents the extreme point that ensures (\ref{constr1}) holds. Consequently, $\alpha_i$ must satisfy $\rho d_{i,j} \leq \alpha_i < s_i$ for all $j \in \mathcal{J}$. It should be noted that the problem in \eqref{dual_problem} is a maximization problem, and $\alpha_i$ reaches its extreme point at a higher value. Therefore, for each area $i$, we have a set of extreme points if there exists an EN $j^* \in \mathcal{J}$ such that $d_{i,j} < d_{i,j^{*}}$ for all $j \in \mathcal{J} \setminus \{j^*\}$, and $\alpha_i = \rho d_{i,j^{*}}$. 

Overall, the optimal objective value for \textbf{Case 2} is given as:
\begin{align}
\label{optimal_obj_case2}
    \rho d_{i,j^{*}} \lambda_i(y) + \sum_{j: d_{i,j} < d_{i,j^{*}} } C_{i,j} \rho ( d_{i,j} -  d_{i,j^{*}}) y_j, ~ \forall i. 
\end{align}

Since $s_i > \rho d_{i,j^{*}}$ and $\rho - \frac{s_i}{d^{\sf min}_{i}} < \rho - \frac{\rho d_{i.j^{*}}}{d^{\sf min}_{i}}$, by considering these two cases together, we can derive a closed-form expression for the optimal objective value of the model in (\ref{dual_problem}). For a given $j^*$, this is achieved by taking the maximum between the following:
\[\begin{cases} 
       \rho d_{i,j^{*}} \lambda_i(y) + \sum_{j: d_{i,j} < d_{i,j^{*}}} C_{i,j} y_j \rho ( d_{i,j} -  d_{i,j^{*}})  \nonumber  \\
      \rho d_{i,j^{*}} \lambda_i(y) + \bigg[ \bigg(\rho - \frac{\rho d_{i,j^{*}}}{d_{i}^{\sf min}}\bigg) \Delta_i \bigg] \lambda_i(y).
\end{cases} \]
As the dual problem (\ref{dual_problem}) is both feasible and bounded, strong duality holds, and the expression given above represents the optimal objective value for the primal problem $g_i(y, \lambda)$.

Furthermore, for each area $i \in \mathcal{I}$, the optimal value of the inner problem $g_i(y, \lambda)$ corresponding to the actual realization $\xi_n \in \Xi$ with the probability $p_{i,n}$ can be given as follows:
\begin{align}
\label{optimal_obj}
    \theta_{i,n} (y) = & \max_{j^{*}\in\mathcal{J}} ~~ \rho d_{i,j^{*}} \xi_n + \max \bigg\{ \big[ \rho - \frac{\rho d_{i,j^{*}}}{d_{i}^{\sf min}}) \Delta_i \big] \xi_n, \nonumber \\
    & \sum_{j: d_{i,j} < d_{i,j^{*}} } C_{i,j} \rho ( d_{i,j} -  d_{i,j^{*}})  y_j \bigg\}, ~\forall i,n.
\end{align}
Both terms within the inner maximization of (\ref{optimal_obj}) are negative. Intuitively,  the inner objective determines which one of these negative terms imposes a more stringent requirement, either in terms of the capacity constraint or the delay constraint.

\textbf{\textit{Step 2:}} \textit{Step 1} reformulates the inner bilevel max-min problem to a single-level maximization. Given the EN placement decision $y$, we can reformulate the inner problem as $\max_{\pi \in \mathcal{U}_2(y)} \mathbb{E} \big[ g(y, \lambda(y)) \big]$, as shown below:
\begin{subequations}
\label{dro_step2}
\begin{align}
& \max_{p_{i,n}} ~~ \sum_{i \in \mathcal{I}} \sum_{n = 1}^{N} p_{i,n}  \theta_{i,n} (y) \\ 
\text{s.t.} ~~~  & \sum_{n = 1}^{N }p_{i,n} \xi_{n} = 1, ~ \forall i ~~~~(\omega_i) \\
& \sum_{n = 1}^{N }p_{i,n} \xi_{n} \leq \Gamma_i^{\mu} + \mu_i (y), ~ \forall i ~~~~(\delta_i^{1})\\
& \sum_{n = 1}^{N }p_{i,n} \xi_{n} \geq \Gamma_i^{\mu} - \mu_i (y), ~ \forall i ~~~~(\delta_i^{2})\\
& \sum_{n = 1}^{N } p_{i,n} \xi_{n}^2 \leq \bigg( \sigma_i^{2} + (\mu_i({y}))^2 \bigg) \bar{\Gamma}^{\sigma}_{i}, ~ \forall i  ~~~~ (\gamma_i^{1})\\
& \sum_{n = 1}^{N } p_{i,n} \xi_{n}^2 \geq \bigg( \sigma_i^{2} + (\mu_i({y}))^2 \bigg) 
\underline{\Gamma}_{i}^{\sigma},  ~ \forall i, ~~~~ (\gamma_i^{2})
\end{align}
\end{subequations}
where ($\omega_i$,$\delta_i^{1}$,$\delta_i^{2}$,$\gamma_i^{1}$,$\gamma_i^{2}$), for all $i$, are the dual variables associated with all constraints. According to duality theory for LP \cite{RObook}, we can formulate the corresponding dual problem as:
\begin{subequations}
\label{dro_step2}
\begin{align}
& \min_{\omega,\delta,\gamma} ~~ \sum_{i} \omega_i + \delta_i^1 (\mu_i(y) + \Gamma_i^{\mu}) - \delta_i^{2} (\mu_i(y) - \Gamma_i^{\mu}) +  \nonumber \\
& \big( \sigma_i^{2}(y) + ((\mu_i({y}))^2 \big) \bar{\Gamma}^{\sigma}_{i} \gamma_i^1 -  \big( \sigma_i^{2} (y) + (\mu_i({y}))^2 \big) \underline{\Gamma}^{\sigma}_{i} \gamma_i^2 \\ 
&\text{s.t.} ~~~ \omega_{i} + (\delta_{i}^1 -  \delta_{i}^2) \xi_n + (\gamma_i^{1} - \gamma_i^{2}) \xi_{n}^{2} \geq  \theta_{i,n} (y) \\
& \qquad \delta_i^{1}, \delta_i^{2}, \gamma_i^1, \gamma_i^2 \geq 0, ~\forall i.
\end{align}
\end{subequations}
Recall that our decision-dependent ambiguity set is constructed based on the mean $\mu_i(y)$ and variance $\sigma_i^2(y)$, which are defined as affine functions of the decision $y$, as shown in (\ref{moment_def}). Hence, $(\mu_{i}(y))^2$ can be extended to the following form:
\begin{align*}
    \bar{\mu}_i^2 \bigg( 1 + \sum_{j \in \mathcal{J}} (2 \Psi_{i,j}^{\mu} + (\Psi_{i,j}^{\mu})^2 ) y_{j} 
    + 2 \sum_{l = 1}^{J} \sum_{m = 1}^{l -1} \Psi_{i,l}^{\mu} \Psi_{i,m}^{\mu} y_m y_l  \bigg).
\end{align*}
By incorporating the definitions of $(\mu_i(y))$, $(\sigma_i(y))$, and $(\mu^2_i(y))$, we obtain the following MINLP:
\begin{subequations}
\label{mnlp}
\begin{align}
& \min_{y \in \{0,1\}}  ~~ \sum_{j} f_j y_{j} + \sum_{i} \bigg( \omega_i + \delta_i^1 (\bar{\mu}_i + \Gamma_{i}^{\mu}) - \delta_i^2 (\bar{\mu}_i + \Gamma_{i}^{\mu}) \nonumber \\
&\qquad + \bar{\mu}_i \sum_{j} \Psi_{i,j}^{\mu} ( \delta_i^1 y_j - \delta_i^2 y_j ) + (\bar{\sigma}_i^2 + \bar{\mu}_{i}^2) (\bar{\Gamma}_i^{\sigma} \gamma_i^1 - \underline{\Gamma}_i^{\sigma} \gamma_i^2 ) \nonumber \\ 
&\qquad + 2 \bar{\mu}_i^2  \sum_{l = 1}^{J} \sum_{m = 1}^{l-1} \Psi_{i,l}^{\sigma} \Psi_{i,m}^{\sigma}   \big( \bar{\Gamma}_i^{\sigma} y_l y_m \gamma_i^1 - \underline{\Gamma}_i^{\sigma} y_l y_m \gamma_i^2  \big) \nonumber\\
&\qquad +  \sum_{j} \Lambda_{i,j} \big( \bar{\Gamma}_i^{\sigma} \gamma_i^1 y_j - \underline{\Gamma}_i^{\sigma}   \gamma_i^2 y_j \big)  \\
& \text{s.t.} ~~~ (\ref{budget}), ~ (\ref{reliability})\\
& \Lambda_{i,j} = \bar{\mu}_i^{2} \big( (\Psi_{i,j}^{\mu})^2  + 2 \Psi_{i,j}^{\sigma}  \big) - \bar{\sigma}_i^{2} \Psi_{i,j}^{\mu}, ~\forall i,j  \\
& \omega_{i} + (\delta_{i}^1 -  \delta_{i}^2) \xi_n + (\gamma_i^{1} - \gamma_i^{2}) \xi_{n}^{2} \geq  \theta_{i,n} (y), ~ \forall i,n \label{feasible_constr1} \\
& \delta_i^{1}, \delta_i^{2}, \gamma_i^1, \gamma_i^2 \geq 0, ~ \forall i . \label{feasible_constr2}
\end{align}
\end{subequations}

\textbf{\textit{Step 3:}} After \textbf{Step 2}, the proposed problem is reformulated into a single-level MINLP, which contains multiple trilinear and bilinear terms. It is important to note that all bilinear and trilinear terms in our formulation exhibit similar characteristics: they involve the product of binary variables and a non-negative continuous variable. To address these nonlinear relationships, we employ McCormick envelopes for linearization. Due to the space limitation, we present the linearization of one bilinear term and one trilinear term. The linearization process for the remaining terms follows a similar approach.

For notation brevity, we use $r \in \{1,2\}$ as the superscript on the dual variable. Let $\mathcal{M}_{\kappa,y,\gamma}$ denote the set involving the McCormick inequalities for linearizing any bilinear term ($\kappa^{r} = \gamma^{r} y$), where $y \in \{0,1\}$, and $\gamma^r$ is non-negative. We have
\begin{align}
\label{bilinear}
    \mathcal{M}_{\kappa, y,\gamma} = \bigg\{& (\kappa, \gamma, y): ~ \underline{\gamma}^{r} y \leq \kappa^{r} \leq \bar{\gamma}^{r} y, \underline{\gamma}^r \leq \gamma^r \leq \bar{\gamma}^r \nonumber \\
    & \gamma^{r} - (1 - y) \bar{\gamma}^{r} \leq \kappa^{r} \leq \gamma^{r} - (1 - y) \underline{\gamma}^{r}  \bigg\}, 
\end{align}
where $\Bar{\gamma}^r$ and $\underline{\gamma}^r$ are the upper bound and lower bound, respectively, on the dual variable $\gamma^r$. In our problem, we set these bounds as sufficiently large positive numbers, denoted by $M$. Similarly, we denote $\mathcal{M}_{\eta,y_l,y_m,\gamma}$ as the set of McCormick inequalities used to linearize trilinear terms. A trilinear term in our formulation contains one non-negative variable and two binary variables. For example, $\eta^{r}_{i,l,m} = \gamma_i^r Y_{l,m} = \gamma_i^r y_l y_m$ has $\gamma_i^r \geq 0$ and $y_l,y_m \in \{0,1\}$. $\mathcal{M}_{\eta,y_l,y_m,\gamma}$ is given as follows:
\begin{align}
\label{trilinear}
    & \mathcal{M}_{\eta,y_l,y_m,\gamma} = \bigg\{ (\eta,\gamma, y_l, y_m): ~\eta^{r} \leq \bar{\gamma}^r y_l, \eta^{r} \leq \bar{\gamma}^r y_m, \nonumber \\
    &\qquad\qquad \eta^{r} \leq \gamma^r - (1 - y_l) \underline{\gamma}^r, ~\eta^{r} \leq \gamma^r - (1 - y_m) \underline{\gamma}^r,  \nonumber  \\
    &\qquad\qquad \eta^{r} \geq \underline{\gamma}^r (y_l + y_m -1), ~\eta^{r} \geq \gamma^r + \bar{\gamma}^r (y_l + y_m - 2), \nonumber\\ 
    &\qquad\qquad y_l \leq 1, ~y_m \leq 1, ~\underline{\gamma}^r \leq \gamma^r \leq \bar{\gamma}^r \bigg\}. 
\end{align}
According to the McCormick linearization for bilinear and trilinear terms, we can derive the following MILP formulation ($\mathcal{P}_1^{'}$), which provides an exact solution for the proposed problem under the decision-dependent ambiguity set (\ref{ddu}):
\begin{subequations}
\label{milp}
\begin{align}
    & (\mathcal{P}_1^{'}):  \min  \sum_{j} f_j y_{j} +\! \sum_{i}\! \Big( \omega_i + \delta_i^1 (\bar{\mu}_i + \Gamma_{i}^{\mu}) - \delta_i^2 (\bar{\mu}_i - \Gamma_{i}^{\mu}) \nonumber \\
    & \!\!\qquad\quad + \bar{\mu}_i \sum_{j} \Psi_{i,j}^{\mu} ( \tau_{i,j}^1 - \tau_{i,j}^2 ) \!+\! (\bar{\sigma}_i^2 + \bar{\mu}_{i}^2) (\bar{\Gamma}_i^{\sigma} \gamma_i^1 - \underline{\Gamma}_i^{\sigma} \gamma_i^2 ) \nonumber \\ 
    & \!\!\qquad\quad+ 2 \bar{\mu}_i^2  \sum_{l = 1}^{J} \sum_{m = 1}^{l-1} \Psi_{i,l}^{\sigma} \Psi_{i,m}^{\sigma}   \big( \bar{\Gamma}_i^{\sigma} \eta_{i,l,m}^1  - \underline{\Gamma}_i^{\sigma} \eta_{i,l,m}^2  \big) \nonumber\\
    & \!\!\qquad\quad+  \sum_{j} \Lambda_{i,j} \big( \bar{\Gamma}_i^{\sigma} \kappa_{i,j}^1 - \underline{\Gamma}_i^{\sigma}   \kappa_{i,j}^2 \big)\Big)  \\
    & \text{s.t.} ~~ (\ref{budget}), ~ (\ref{reliability})  \\
    &\qquad  \Lambda_{i,j} = - \bar{\sigma}_i^{2} \Psi_{i,j}^{\mu} + \bar{\mu}_i^{2} \big( (\Psi_{i,j}^{\mu})^2  + 2 \Psi_{i,j}^{\sigma}  \big), ~\forall i,j  \\
    &\qquad \omega_{i} + (\delta_{i}^1 -  \delta_{i}^2) \xi_n + (\gamma_i^{1} - \gamma_i^{2}) \xi_{n}^{2} \geq  \theta_{i,n} (y), ~\forall i, n \\
    &\qquad \kappa^{r}_{i,j} \in \mathcal{M}_{y_j,\gamma^{r}_i}, ~ \tau^{r}_{i,j} \in \mathcal{M}_{y_j,\gamma^{r}_i}, ~\forall i,j,r \\
    &\qquad \eta^{r}_{i,l.m} \in \mathcal{M}_{y_l,y_m,\gamma^{r}_i}, ~\forall i,l,m,r, l > m \\ 
    &\qquad y_j \in \{0,1\}, ~\forall j; ~ \delta_i^{1}, \delta^{2}, \gamma_i^1, \gamma_i^2 \geq 0, ~\forall i. 
\end{align}
\end{subequations}
Hence, by following \textbf{\textit{Step 1}} to \textbf{\textit{Step 3}}, we can achieve an exact MILP reformulation that can be solved directly using off-the-shelf solvers. These steps are summarized in \textbf{Algorithm \ref{ALG:OPT_placement}}.

\begin{algorithm}[ht!]
\caption{Exact OPT-Placement}
\label{ALG:OPT_placement}
\begin{algorithmic}[1]
\STATE \textbf{Initialization:} $f_{j}$, $s_{i}$, $d_{i,j}$, $\mathcal{U}_2(y)$, $B$, $K^{\sf min}$, $\rho$
\STATE \textbf{Step 1:} Solve the inner problem (\ref{primal_1}): Find analytical expression (\ref{optimal_obj}) for the dual problem of inner objective (\ref{dual_problem}). 
\STATE \textbf{Step 2:} Dualize the obtained inner ``\textit{max-max}" problem which is subject to all constraints within $\mathcal{U}_2(y)$. 
\STATE \textbf{Step ~3:} Linearize the bilinear and trilinear terms based on McCormick linearization techniques.
\STATE \textbf{Output:} Optimal placement $(y^*)$.
\end{algorithmic}
\end{algorithm}

\subsection{Improved Variant}
While \textbf{Algorithm~\ref{ALG:OPT_placement}} provides an optimal solution to the proposed problem ($\mathcal{P}_1$), the computational time can become sensitive to the network size. To address this limitation, the improved algorithm takes a different approach. It generates extreme rays for the feasible region of ($\mathcal{P}_1$), which are then incorporated into the MILP problem obtained in ($\mathcal{P}_1^{'}$). By identifying the potential locations of extreme points for these dual variables, we achieve a stronger reformulation compared to ($\mathcal{P}_1^{'}$). This enhancement helps reduce computational complexity and allows for more efficient solving of larger networks.

After completing \textbf{Step 2}, formulation (\ref{mnlp}) shows that the problem is feasible within a region satisfying the inequalities (\ref{feasible_constr1}) - (\ref{feasible_constr2}). To this end, let us define $\delta_i = \delta_i^1 - \delta_i^2$ and $\gamma_i = \gamma_i^1 - \gamma_i^2$. It should be noted that $\delta_i$ and $\gamma_i$ are unbounded, which means searching for the extreme point to achieve the optimum objective might be time-consuming. Therefore, the objective is to determine a set of extreme rays ($\omega_i$, $\delta_i^1$, $\delta_i^2$, $\gamma_i^1$, $\gamma_i^2$) that effectively represent the feasible region defined by constraints (\ref{feasible_constr1}) and (\ref{feasible_constr2}). As a result, the improved variant can strengthen the problem by reducing the feasible region of the inner problem, leading to faster computation times.

To identify extreme rays, we solve the following inequality system for $k,l \in \{1,2,\dots, N \}$, where $k$ and $l$ represent the indices of extreme points: 
\begin{subequations}
\begin{align}
    & \omega_i + \delta_i \xi_k + \gamma_j \xi_k^2 = 0, ~ \forall i,k \label{linear_1}\\ 
    & \omega_i + \delta_i \xi_l + \gamma_j \xi_l^2 = 0, ~\forall i,l \label{linear_2}\\ 
    & \omega_i + \delta_i \xi_n + \gamma_j \xi_n^2 \geq 0, ~\forall n \in \{1,2,\dots,N\} \setminus \{l,k\}.  \label{linear_3}
\end{align}
\end{subequations}
Without loss of generality, we assume that $\xi_k < \xi_l$. Our objective is to determine the relationship between $\xi_k$, $\xi_l$, and the other instances $\xi_n, n \in \{1,2,\dots, N\} \setminus \{k,l\}$. To achieve this, we define $\{ \xi_{(1)}, \xi_{(2)}, \dots, \xi_{(N)} \}$ as a ordered support for the random demand.

According to (\ref{linear_1}) and (\ref{linear_2}), we can derive $\delta_i = -(\xi_k + \xi_l) \gamma_i$ and $\omega_i = \xi_k \xi_l \gamma_i$. To ensure clarity, we will fix the direction of the unit vector along $\gamma_{i}$ and determine the direction of other variables to satisfy the inequality (\ref{linear_3}). By normalizing $\gamma_i$, we have $\abs{\gamma_i} = 1$, which will be analyzed in the following: 

\noindent
\textbf{Case 1 ($\gamma_ i = 1$):} Due to the assumption $\xi_k < \xi_l$ and $\xi_n \geq 0$, we must ensure $(\xi_n - \xi_l)(\xi_n - \xi_k) \geq 0, \forall n \in \{1,2,\dots,N \} \setminus \{ k,l \}$. Thus, we can have either $\xi_n \geq \xi_l > \xi_k $ or $\xi_n \leq \xi_k < \xi_l$. 
Based on these two relationships, we can derive the expressions for $\omega_i$ and $\delta_i$. 
Thus, there exists two extreme rays ($\omega_i$,$\delta_i$,$\gamma_i$) that satisfy the conditions: 
\begin{subequations}
\begin{align}
    & \omega_i =  \xi_{(1)} \xi_{(2)}, ~\delta_i = -(\xi_{(1)} + \xi_{(2)}), ~\gamma_i = 1, \forall i\\
    & \omega_i = \xi_{(N-1)} \xi_{(N)}, ~\delta_i = -(\xi_{(N-1)} + \xi_{(N)}), ~\gamma_i = 1, \forall i.
\end{align}
\end{subequations}
\noindent
\textbf{Case 2 ($\gamma_ i = -1$):} Due to the assumption $\xi_k < \xi_l$ and $\xi_n \geq 0$, we must ensure $(\xi_n - \xi_l)(\xi_n - \xi_k) \leq 0, \forall n \in \{1,2,\dots,N \} \setminus \{ k,l \}$. Thus, we have  $\xi_k \leq \xi_n \leq \xi_l $. Therefore, the extreme ray can be expressed as: 
\begin{align}
    \omega_i =  \xi_{(1)} \xi_{(N)}, ~\delta_i = -(\xi_{(1)} + \xi_{(N)}), ~\gamma_i = -1, ~ \forall i.
\end{align}
Given the values of $\gamma_i$ and $\delta_i$, we can express $\delta_i^1 = \max \{0,\delta_i\}$, $\delta_i^{2} = \max \{0,-\delta_i\}$, $\gamma_i^1 = \max \{0,\gamma_i\}$ , and $\gamma_i^2 = \max \{0,-\gamma_i\}$. Thus, we can substitute ($\omega_i$,$\delta_i^1$,$\delta_i^2$,$\gamma_i^1$,$\gamma_i^2$) into the problem (\ref{mnlp}).
The following inequalities ensure that the dual problem (\ref{dro_step2}) is bounded, thereby guaranteeing the feasibility of ($\mathcal{P}_1$). Thus, for every area $i \in \mathcal{I}$, we have:
\begin{subequations}
\label{valid_ineq}
\begin{align}
    & \!\!\xi_{(1)} \xi_{(2)} - (\xi_{(1)} + \xi_{(2)}) (\mu_i(y) -  \Gamma_{i}^{\mu}) + S_i(y) \bar{\Gamma}_i^{\sigma}  \geq  0 \\ 
    & \!\!\xi_{(N-1)} \xi_{(N)} \!-\! (\xi_{(N-1)} \!+\! \xi_{(N)}) (\mu_i(y) \!-\! \Gamma_{i}^{\mu}) \!+\! S_i(y) \bar{\Gamma}_i^{\sigma} \!\geq \!0\!\!\\ 
    &  \!\!- \xi_{(1)} \xi_{(N)} + (\xi_{(1)} + \xi_{(N)}) (\mu_i(y) -  \Gamma_{i}^{\mu}) \!-\! S_i(y) \underline{\Gamma}_i^{\sigma} \geq 0, \!\!
\end{align}
\end{subequations}
where $S_i(y) = \sigma_i^2(y) + \mu^2_i(y)$. Similar to the previous section, we employ McCormick linearization techniques to linearize the bilinear terms in (\ref{valid_ineq}). As a result, (\ref{valid_ineq}) becomes:
\begin{subequations}
\label{valid_ineq_linear}
\begin{align}
    & \!\!\xi_{(1)} \xi_{(2)} - (\xi_{(1)} + \xi_{(2)}) (\zeta_i - \Gamma_i^{\mu}) + z_i \bar{\Gamma}_{i}^{\sigma} \geq 0, ~ \forall i\\
    & \!\!\xi_{(N-1)} \xi_{(N)} \!-\! (\xi_{(N-1)} \!+\! \xi_{(N)})  (\zeta_i - \Gamma_i^{\mu}) \!+\! z_i \bar{\Gamma}_i^{\sigma} \geq 0, ~\forall i\\ 
    &  \!\!- \xi_{(1)} \xi_{(N)} + (\xi_{(1)} + \xi_{(N)}) ( \zeta_i + \Gamma_{i}^{\mu}) - z_i \underline{\Gamma}_i^{\sigma} \geq 0, ~ \forall i \\
    &\!\! z_i \!=\!  \Bar{\sigma}_i^2  \!+\! \Bar{\mu}_i^2 \!+\!\! \sum_{j} \Lambda_{i,j} y_j \!+\! 2 \Bar{\mu}_i^2 \!\sum_{l = 1}^{J} \!\sum_{m = 1}^{l-1} \!\Psi_{i,l}^{\mu} \Psi_{i,m}^{\mu} Y_{l,m}, ~\forall i\\
    &\!\! \zeta_i = \bar{\mu}_i (1 + \sum_{j} \Psi_{i,j}^{\mu} y_j), ~ \forall i \\
    &\!\!  Y_{l,m} \geq  y_{l} + y_{m} - 1, Y_{l,m} \leq y_l, Y_{l,m} \leq y_m, ~\forall l,m.
\end{align}
\end{subequations}
After incorporating the constraints in (\ref{valid_ineq_linear}) into the MILP reformulation ($\mathcal{P}_1^*$), we strengthen our formulation for the proposed problem by reducing the feasible region. This enhancement allows for efficient computation and a tighter representation of feasible solutions.

\vspace{-0.2cm}
\section{Numerical Results}
\label{results}
\subsection{Simulation Setting}
\label{simsetting}
We consider an EC system with $I = 15$ areas and $J = 10$ ENs in the default setting while larger networks will also be considered in sensitivity analyses. The edge network topology is generated based on the cities and locations of randomly selected Equinix edge data centers (DCs) \cite{equinix}. The network delay ($d_{i,j}$) between any two selected DCs is obtained directly from the global ping dataset \cite{wondernetwork}. The EN placement cost $f_j$ is sampled from the uniform distribution $U(20,30)$, while the unmet penalty $s_i$ is randomly generated from $U(30,40)$. To generate $C_{i,j}$, the maximum resource capacity ($C^{\sf max}_j$) at each EN $j$ is randomly selected from the set $\{84,96,128\}$ \textit{vCPUs}. The resources available to each area $i$, $C_{i,j}$, are pre-allocated according to the relative historical demand. 

We randomly generate the empirical mean of resource demand ($\bar{\mu}_i$) in each area $i$, following a uniform distribution with values ranging from $20$ to $50$ \textit{vCPUs}. We define $\theta_i$ as the ratio of variation at each area, represented by $\theta_i = \frac{\bar{\sigma}_i}{\bar{\mu}_i}, \forall i$. The support size of demand ($N$) at each area is taken as $100$, with $\xi_1, \dots, \xi_N$ in the range $\{1,\dots,100\}$. To establish decision dependency between the demand distribution and EN placement decisions, $\Psi_{i,j}^{\mu}$ and $\Psi_{i,j}^{\sigma}$ are considered as decreasing functions of the corresponding network delay (e.g., distance), i.e., $\exp\left(-\frac{ d_{i,j}}{b}\right), \forall i,j$, where $b$ is a parameter controlling the decaying rate. This means the placement of EN $j$ has a higher impact when it is closer to nearby areas. Since $\sum_{j} \Psi_{i,j}^{\sigma} \leq 1$, we normalize both impact parameters. To control the level of robustness with respect to the true mean and true variance of demand, we define $\mathbf{\epsilon}_i = (\epsilon_i^{\mu}, \epsilon_i^{\sigma}) \in [0,1]^{I}$ to adjust the robustness for each area, i.e., $\Gamma_{i}^{\mu} = \epsilon^{\mu} \mu_i(y)$, $\underline{\Gamma}_i^{\sigma} = 1 - \epsilon_i^{\sigma}$, $\bar{\Gamma}_i^{\sigma} = 1 + \epsilon_i^{\sigma}$. 
In our \textbf{default setting}, the  other system parameters are: $\rho = 0.001$, $\Delta_i = \Delta = 35$, $B = 100$, $b = 25$, $\epsilon_i^{\mu} = \epsilon^{\mu}= 0.8$, $\underline{\Gamma}_i^{\sigma} = \underline{\Gamma}^{\sigma} = 0.8$, $\bar{\Gamma}_i^{\sigma} = \bar{\Gamma}^{\sigma} = 1.2$, $K^{\sf min} = 1$, $\theta_i = \theta = 0.5, \forall i$. We will also vary these important parameters during sensitivity analyses. All the experiments are implemented in MATLAB  using CVX \cite{CVX} and Gurobi \cite{Gurobi}  on a desktop with an Intel Core i7-11700KF and 32 GB of RAM.
\vspace{-0.7cm}
\subsection{Sensitivity analysis}
\label{sensi}
This section conducts sensitivity analyses to assess the impact of key system parameters on the optimal solution. The parameters under investigation include the budget ($B$), delay penalty ($\rho$), and impact factors $\Psi_{i,j} = (\Psi_{i,j}^{\mu},\Psi_{i,j}^{\sigma})$. To evaluate the impact of the EN placement cost $\bf{f}$, we introduce a scaling factor $h$, where $h = 1$ represents the default setting. The base value of $\bf{f}$ generated in Section \ref{simsetting} is multiplied by $h$ to either scale up or down the placement cost. A higher value of $h$ indicates a higher EN placement cost.\vspace{-0.4cm}
\begin{figure}[h!] 
\centering
         \subfigure[Cost: varying $h$ and $B$]{
	    \includegraphics[width=0.242\textwidth,height=0.102\textheight]{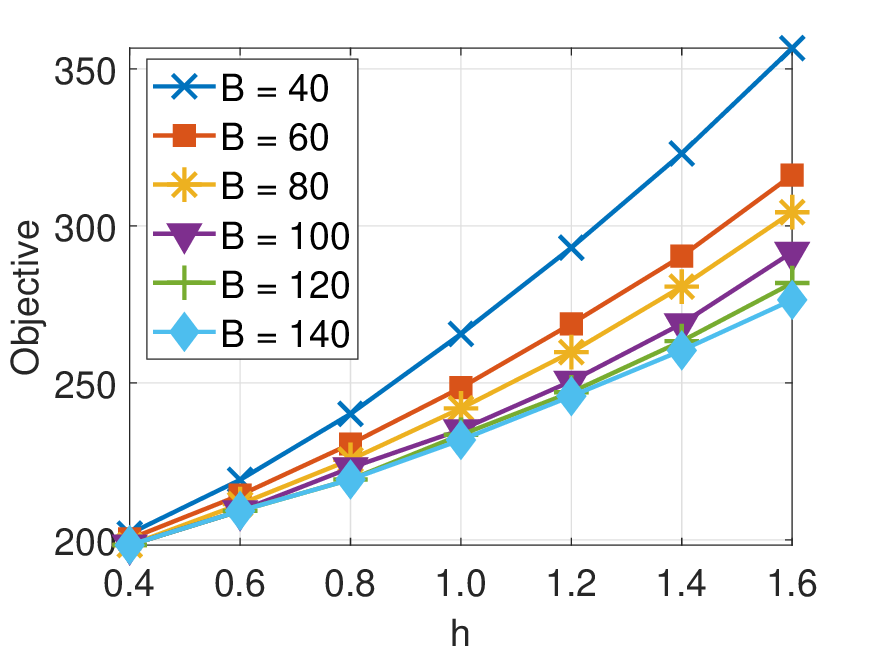}
	     \label{fig:B_h_cost}
	}  \hspace*{-2.1em} 
	     \subfigure[Payment: varying $h$ and $B$]{
	     \includegraphics[width=0.242\textwidth,height=0.102\textheight]{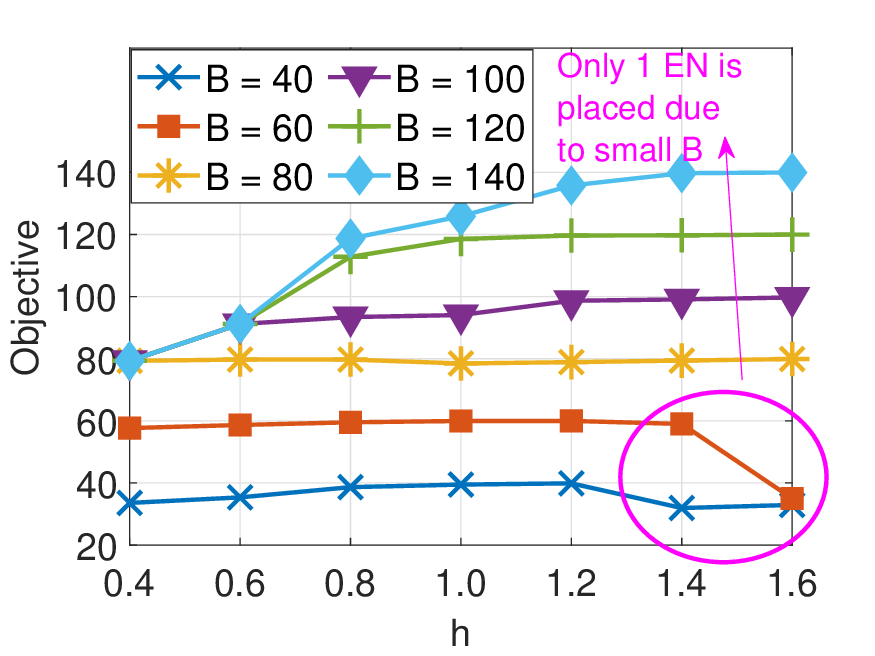}
	     \label{fig:B_h_payment}
	}
         \subfigure[Varying $\Delta$ and $\rho$]{
	    \includegraphics[width=0.242\textwidth,height=0.102\textheight]{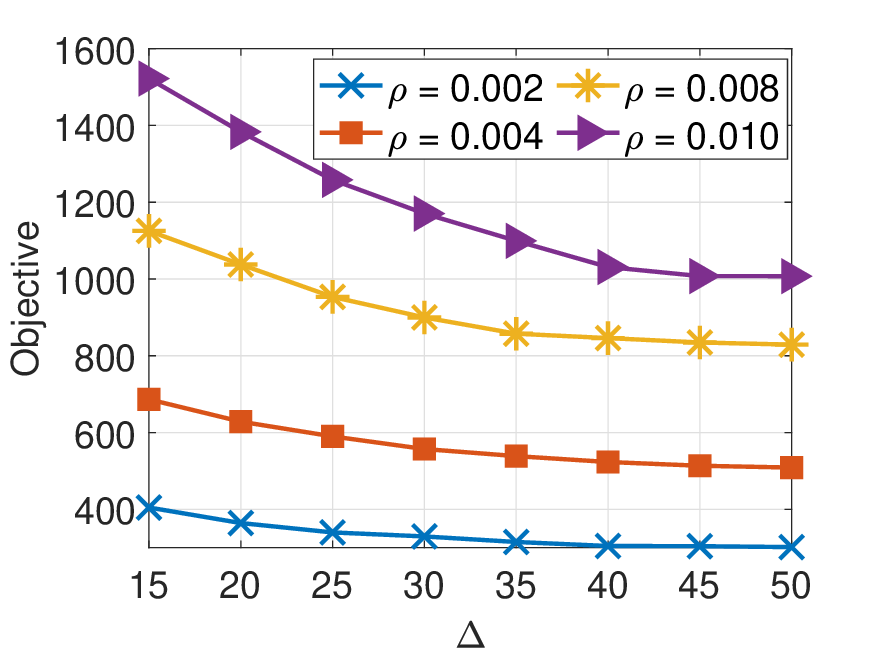}
	     \label{fig:rho_Delta_cost}
	}   \hspace*{-2.1em} 
	    \subfigure[Varying $\epsilon^{\mu}$ and $\epsilon^{\sigma}$]{
	    \includegraphics[width=0.242\textwidth,height=0.102\textheight]{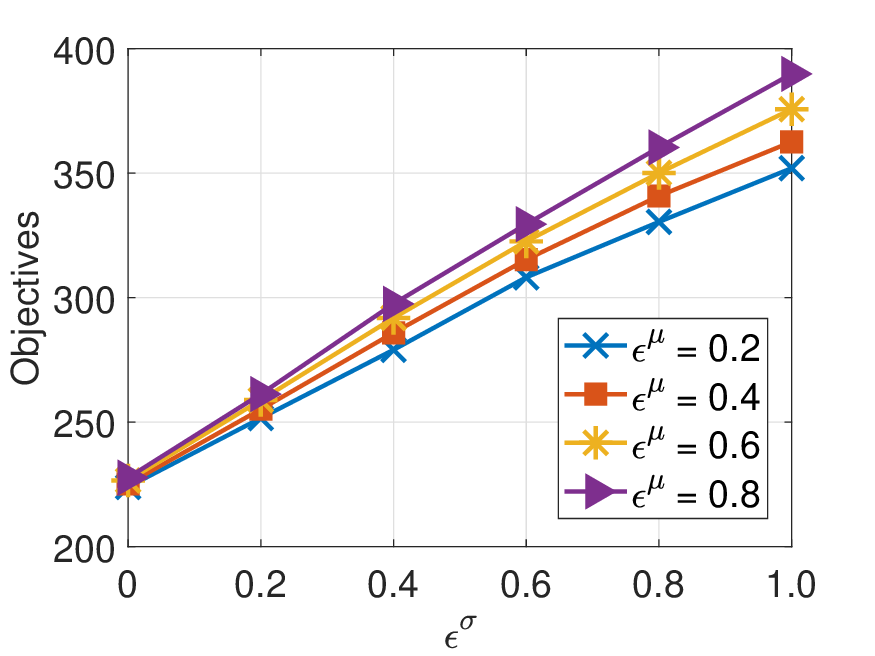}
	     \label{fig:sigma_mu_cost}
	} 
	\vspace{-0.2cm}
	    \caption{Impacts of important system parameters}
     	\vspace{-0.35cm}
\end{figure}

\noindent \textit{\textbf{1) Impacts of placement cost and budget:}} Fig.\ref{fig:B_h_cost} shows that the total cost increases with higher EN placement costs. Due to the limited budget, only a few ENs can be installed, leading the platform to prioritize critical locations for installing ENs, particularly with a higher $h$. However, this may result in increased penalties for unmet demand when $h$ is high. The curve's growth rate increases as $B$ decreases, while increasing $B$ allows the platform to select EN locations more flexibly. Note that the minimum number of placed ENs must be greater than $1$ ($K^{\sf min} = 1$) to avoid situations where the payment equals $0$. Notably, the payment does not monotonically increase with $h$, as highlighted in the dashed circle in Fig.\ref{fig:B_h_payment}. This is because the budget $B = 40$ allows only one EN to be placed with $h = 1.6$, whereas it permits more than 2 ENs with smaller $h$.

\noindent \textit{\textbf{2) Impacts of delay penalty $\rho$ and threshold $\Delta$: }} Fig.~\ref{fig:rho_Delta_cost} demonstrates that the total cost decreases with smaller $\Delta$ and $\rho$ thanks to the improved proximity and reduced delay penalty. The cost reduction helps to mitigate delay and unmet demand penalties, enhancing service quality. The curves converge to a small value after a certain delay threshold, indicating that the delay requirement can be easily satisfied. Thus, the platform begins to prioritize capacity constraints over delay constraints. 

\noindent \textit{\textbf{3) Impacts of parameters in ambiguity set:}} Fig.~\ref{fig:sigma_mu_cost} shows how the level of robustness in the ambiguity set affects the optimal solution. $\epsilon$ controls the distance to the true mean and variance of demand. As $(\epsilon^{\mu}, \epsilon^{\sigma})$ increases, the ambiguity set grows, allowing actual demand to vary over a larger range and incurring a more conservative solution.

\noindent \textit{\textbf{4) Running time comparison:}} We compared the average running time over $20$ randomly generated problem instances for each problem size, as shown in Table~\ref{tab: time_table}. \textbf{Algorithm 1} (denoted as \textbf{\textit{Standard}}) and Algorithm 1 with feasibility cuts (denoted as \textbf{\textit{Improved}}) both produce optimal solutions within a reasonable time for small- and medium-sized networks. However, as the system size grows larger, the advantage of \textit{Improved} over \textit{Standard} becomes more evident. The feasibility cuts reduce the feasible region and mitigate the effects of unbounded dual variables, thus strengthening \textbf{Exact OPT-Placement}. Importantly, it is worth noting that \textit{the underlying problem is a robust planning problem that does not require real-time computation.}

\vspace{-0.2cm}
\begin{table}[h!]
\centering
\begin{tabular}{|c|c|c|}
\hline
\textbf{Network size}  & \textbf{Standard}   & \textbf{Improved}   \\ \hline
I = 10; J = 10 & 31.31s &  21.68s \\ \hline
I = 20; J = 10 & 66.88s &   62.95s\\ \hline 
I = 20; J = 20 &  404.11s & 331.79s   \\ \hline 
I = 30; J = 20 &  1314.8s &  901.8s  \\ \hline
I = 40; J = 20 &  3357.2s &  2178.28s  \\ \hline
\end{tabular}
\caption{Runtime comparison}
\label{tab: time_table}
\vspace{-0.8cm}
\end{table}

\subsection{Performance Comparison}
\label{percom}
In this section, we compare the performance of the proposed \textbf{DRO-DDU} with the following benchmarks:
\begin{itemize}
    \item \textit{HEU}: Choose a subset of ENs according to demand, giving priority to areas with higher demand until the available budget is fully utilized.
    \item \textit{BSPA}: Deploy as many ENs as possible within the budget.
    \item \textit{DET}: Deterministic EN placement problem in (\ref{DET}).
    \item \textit{SO}: Two-stage SO with uniform in-sample distribution. 
    \item \textit{DRO-DIU}: $\Psi_{i,j}^{\mu} \!=\! \Psi_{i,j}^{\sigma} \!=\! 0$. The original problem reduces to a two-stage DRO with exogenous stochastic demand.
\end{itemize}
The platform aims to optimize the EN placement decision ($\by$) before knowing the actual demand. To evaluate the placement decisions provided by different schemes, we conduct an out-of-sample test using model (\ref{Actual_model}). For each scheme, we obtain a placement solution $\hat{\by}$ during the planning stage. Given $\hat{\by}$ and the actual uncertainties $\lambda(\hat{\by})$, the platform can re-optimize the workload allocation decision to minimize the total actual cost. In our experiment, we generate $1000$ scenarios for each scheme to model the actual mean and variance of demand ($\hat{\lambda_i} (\hat{\by})$) satisfying the conditions in (\ref{moment_def}). For each scenario, the platform solves the following actual workload allocation problem (LP):
\begin{subequations}
\label{Actual_model}
\begin{align}
    (\textbf{A}) ~~ \min_{\bx \in \mathbb{R}_{+}^{I \times J} ,\bu \in \mathbb{R}_{+}^{I}}  \quad & \rho \sum_{i,j} d_{i,j} x_{i,j} + \sum_{i} s_i u_i\\
     \text{s.t.} ~~ (\ref{supply_demand}) - (\ref{var_constr1})
\end{align} 
\end{subequations}
The actual total cost is the sum of the EN placement cost and actual workload allocation cost, expressed as:
\begin{align}
\label{Actual_cost}
    \mathcal{C}^{\sf a} = \sum_{j} f_j \hat{y}_{j}  + \rho \sum_{i,j} d_{i,j} x_{i,j}^{\sf a} + \sum_{i} s_i u_{i}^{\sf a},
\end{align}
where ($x^a, u^a$) is the optimal solution to problem (\ref{Actual_model}). The six schemes are evaluated and compared based on their average and worst actual costs over the generated scenarios.

\begin{figure}[h!]
\centering
        \subfigure[Average: varying $\theta$]{
	     \includegraphics[width=0.245\textwidth,height=0.102\textheight]{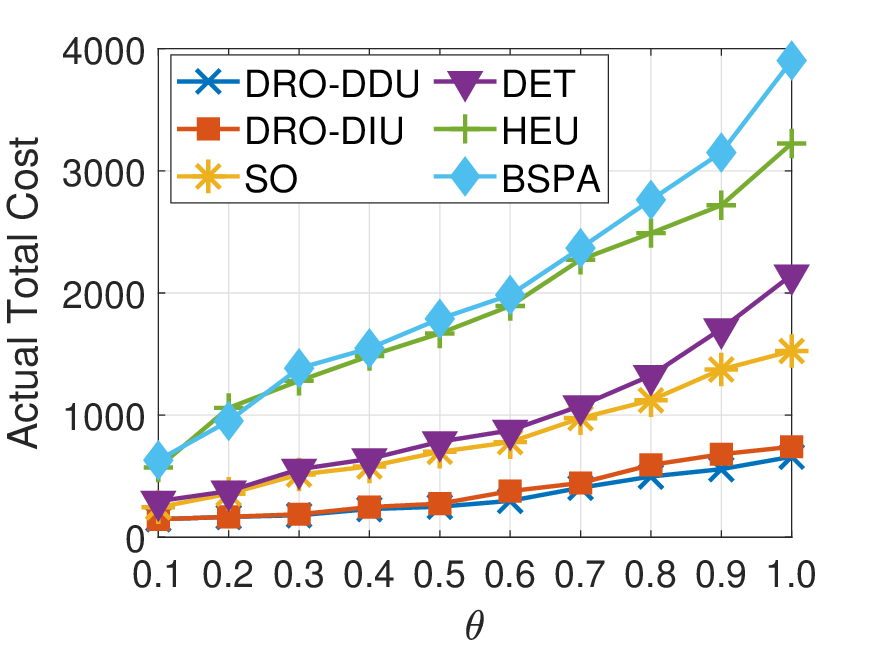}
	     \label{fig:theta_model_avg}
	}   \hspace*{-2.3em}
		\subfigure[Worst: varying $\theta$]{
	     \includegraphics[width=0.245\textwidth,height=0.102\textheight]{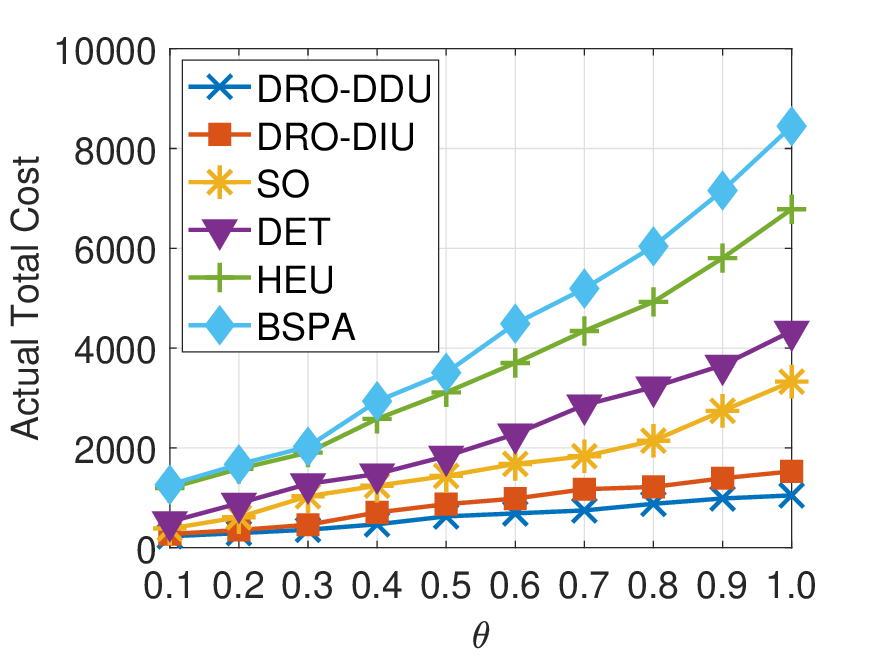}
	     \label{fig:theta_model_worst}
	} 
		 \subfigure[Average: varying $h$]{
	     \includegraphics[width=0.245\textwidth,height=0.102\textheight]{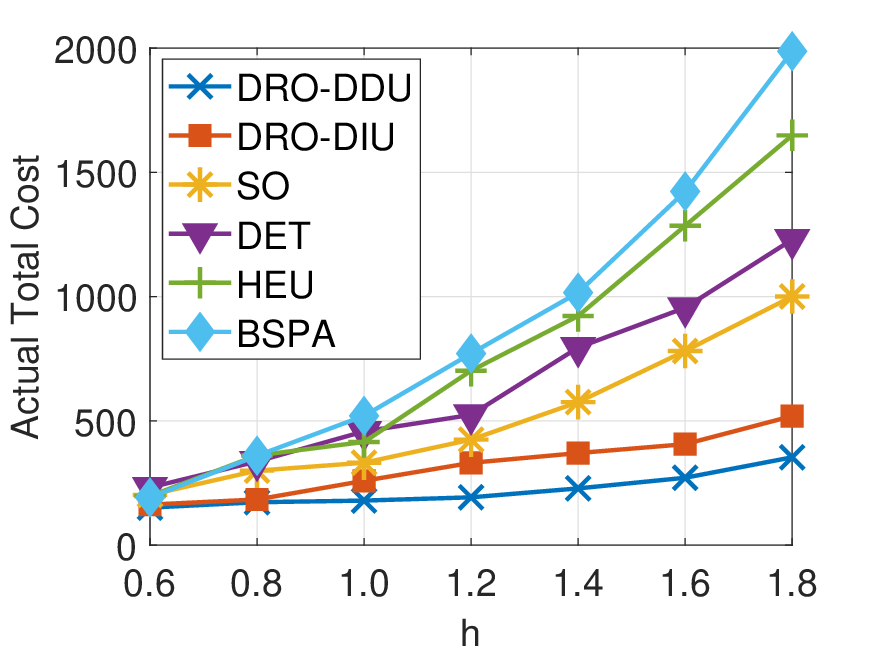}
	     \label{fig:h_cost_avg}
	}   \hspace*{-2.3em}
		\subfigure[Worst: varying $h$]{
	     \includegraphics[width=0.245\textwidth,height=0.102\textheight]{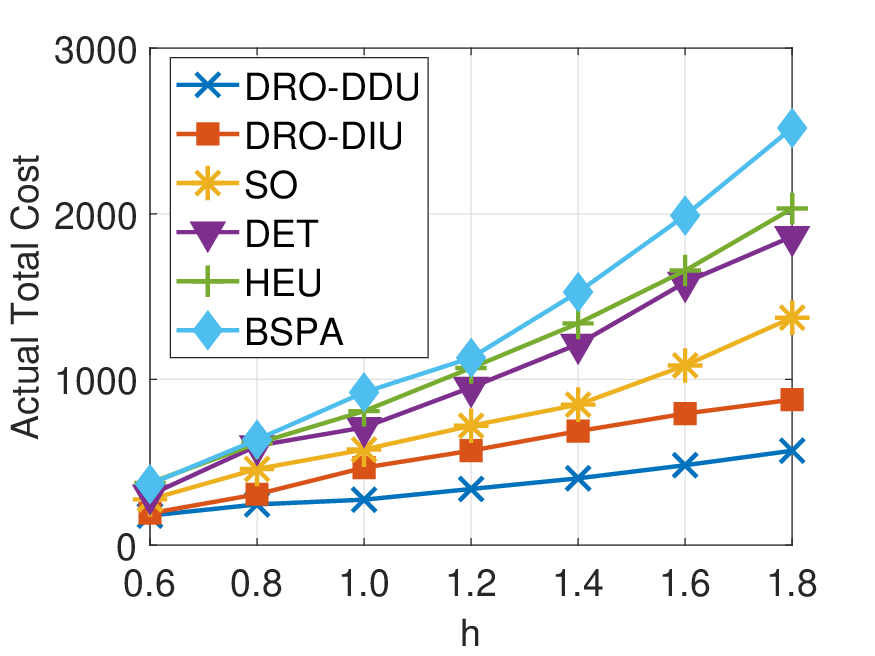}
	     \label{fig:h_cost_worst}
	} 	 	  \vspace{-0.37cm}
	\caption{Model comparison}
	  \vspace{-0.37cm}
\end{figure}

\noindent \textit{\textbf{1) Varying variability:}}
Recall that $\theta_i$ denotes the variation ratio in each area. In Figs.~\ref{fig:theta_model_avg}-\ref{fig:theta_model_worst}, the DRO-based models demonstrate increased stability compared to other schemes, especially with higher variability $\theta_i$. As $\theta$ increases, the gap between these schemes widens due to the significant deviation of actual demand from its mean. \textit{BSPA}, \textit{HEU}, and \textit{DET} do not consider demand uncertainty. \textit{BSPA} performs well with smaller $\theta$ but incurs high costs with a higher $\theta$ due to overly optimistic EN placement. However, the selection of ENs may not be on critical locations, especially within the limited budget. \textit{SO} is prone to out-of-sample disappointment as it relies heavily on in-sample distribution accuracy. In contrast, DRO models consider the worst-case distribution, ensuring robustness in handling various demand patterns during out-of-sample scenarios, as seen in Fig.\ref{fig:theta_model_worst}. \textit{DRO-DDU} model performs significantly better in all settings, highlighting the importance of incorporating decision dependency in uncertainty quantification.
\vspace{-0.38cm}

\noindent \textit{\textbf{2) Varying EN placement cost:}}
Figs.~\ref{fig:h_cost_avg}--\ref{fig:h_cost_worst} demonstrate the superiority of the proposed \textit{DRO-DDU} scheme over other schemes, particularly as the EN placement cost increases. This advantage arises from the proposed model's consideration of decision dependency between planning decisions and uncertain demand, offering better performance compared to alternative schemes. As the cost parameter increases, the performance of solutions obtained from other models monotonically deteriorates. This indicates the critical importance of thoughtfully selecting ENs for planning decisions, particularly when the price of EN placement becomes higher. It can be observed that \textit{BSPA} exhibits decent performance for smaller $h$. This is attributed to the lower $h$, which allows the platform to be relatively optimistic, resulting in the utilization of all budgets for EN placement without careful consideration of critical infrastructures. However, this blind selection worsens its performance as $h$ increases. Additionally, \textit{DRO-DIU} outperforms \textit{SO} since SO relies on the accurate knowledge of the distribution of the historical demand, while the DRO-based model proves more robust in handling various demand patterns.

\begin{figure}[h!]
		 \subfigure[Varying $B$]{
   \includegraphics[width=0.245\textwidth,height=0.102\textheight]{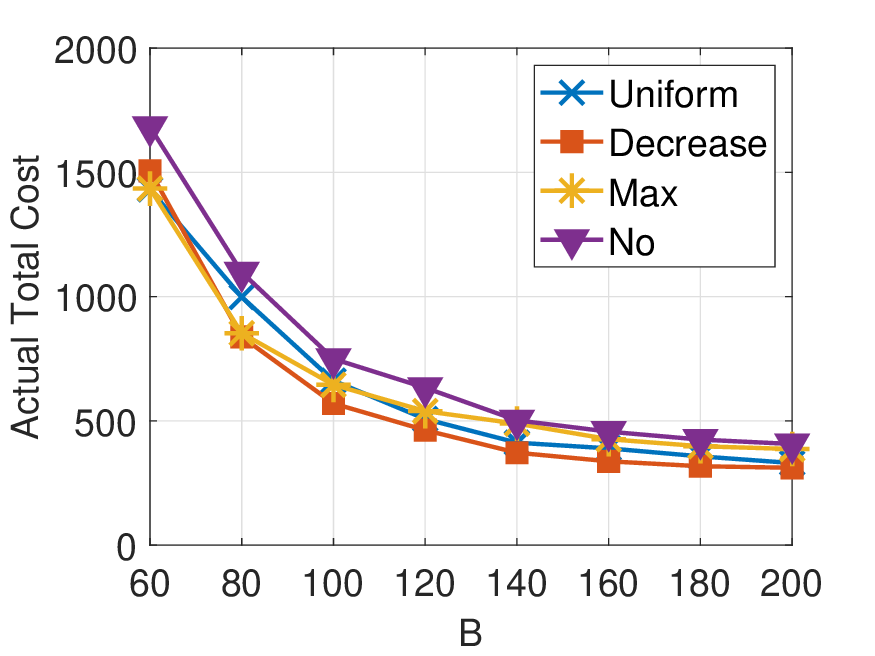}
	     \label{fig:B_DDU_type_cost}
	}   \hspace*{-2.1em} 
	    \subfigure[Varying $b$]{
	    \includegraphics[width=0.245\textwidth,height=0.102\textheight]{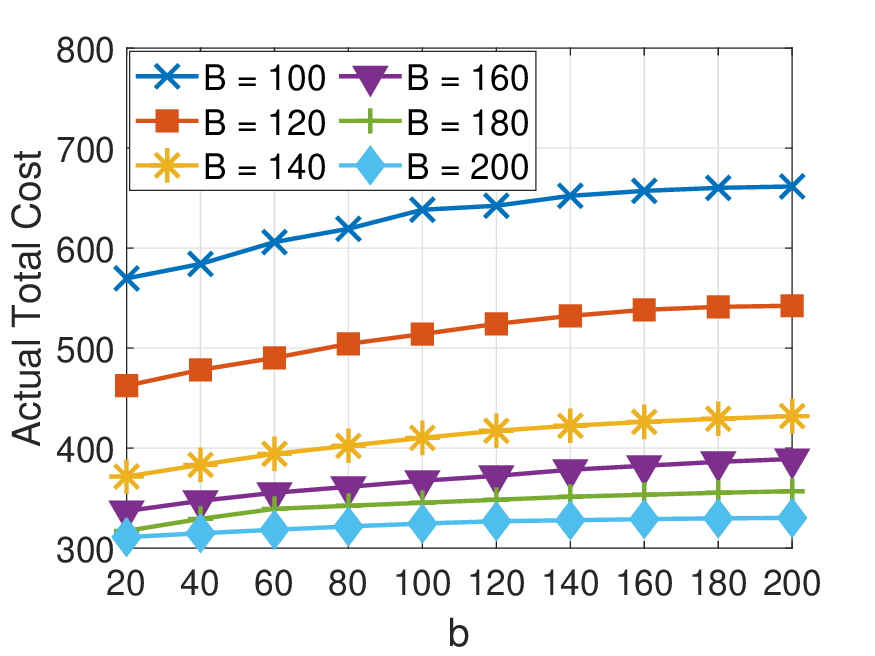}
	     \label{fig:vary_b_cost}
	}      	\vspace{-0.3cm}
	    \caption{Performance comparison with varying DDU sets}
     \vspace{-0.5cm}
\end{figure}

\noindent \textit{\textbf{3) Choice of decision-dependency: }}
Recall that both impact factors ($\mathbf{\Psi}_{i,j}$) are generated from a decreasing function of the network delay between area $i$ and EN $j$, i.e., $e^{-\frac{d_{i,j}}{b}}, \forall i,j$. We denote the proposed form of ambiguity set by \textbf{Decrease}. Additionally, we consider other forms commonly used in practice: (i) (\textbf{Uni}) uniform impact overall areas ($\Psi_{i,j} = \frac{1}{J}$); (ii) (\textbf{No}) No impact: reducing the problem to the traditional DRO problem with a decision-independent ambiguity set; (iii) (\textbf{Max}): Maximum impact on the closest area only ($\min_i d_{i,j}$). Fig.\ref{fig:B_DDU_type_cost} illustrates the impact of decision-dependency in an ambiguity set, influencing the optimal solution with varying $B$. Notably, \textbf{NO} neglects the impact of EN placement decisions, resulting in the highest out-of-sample cost among all choices. On the other hand, \textbf{Max} performs well with the limited budget (i.e., fewer available placed ENs), observing the most significant impact in its nearest area with the highest traffic flow. However, as $B$ becomes larger, the cost exceeds that of \textbf{Uni} and \textbf{Decrease} due to disregarding the impacts of EN placement on demand changes in other areas. \textbf{Uni} assumes uniform impact across all areas but fails to reflect network delays or geographical locations between ENs and nearby areas, resulting in poorer performance in most cases. The proposed impact factor is defined as a decreasing function of the network delay, where $b$ is the decaying rate of the impact for EN $j$ on all areas. A higher value of $b$ corresponds to a slower decaying rate, indicating uniform impact across all areas. Conversely, a smaller value of $b$ leads to a faster-decaying rate, implying that EN $j$ only impacts its closest area. Thus, the proposed form acts as an intermediate state between \textbf{Max} and \textbf{Uni}, and the platform can adjust $b$ to control impact factors according to its budget. Fig.\ref{fig:vary_b_cost} shows that the cost increases with higher $b$ since the impact factor tends to become \textbf{Uni}, failing to reflect the network delays. However, note that the total cost does not always increase monotonically with an increasing $b$. An extremely lower $b$ leads to impact factors similar to \textbf{Max}, performing well only under specific conditions, as discussed previously.

\vspace{-0.4cm}

\section{Conclusion}
\label{conc}
This paper presented a novel two-stage DRO for optimal EN placement, aiming to mitigate the impact of demand uncertainty on user experience. The main novelty of lies in the integration of a DDU set into the DRO framework to capture the interdependence between demand uncertainty and EN placement decisions. To compute the exact optimal solution, two efficient algorithms were developed. Numerical results illustrate the importance and advantages of incorporating endogenous uncertainties, highlighting the improved performance of the proposed model over benchmark schemes.

\bibliographystyle{IEEEtran}
\bibliography{ref.bib}

\appendix

\subsection{Budget-spending priority algorithm (BSPA)}
\label{BSPA_ALG}
This section introduces a \textit{budget-spending priority} benchmark where the platform maximizes the use of these limited funds for EN placement. This scheme is equivalent to solving the following optimization problem: (\ref{BSPA}):
\begin{subequations}
\label{BSPA}
\begin{align}
    & \max_{\bm{y}} ~~ \sum_{j} f_{j} y_{j} \\
    & \textit{s.t} ~~ \sum_{j} f_{j} y_{j} \leq B.
\end{align}
\end{subequations}
The obtained EN placement strategy ($y$) under this scheme will serve as an input to the actual workload allocation model in (\ref{Actual_model}). The actual total operation cost is counted as the sum of the actual workload allocation cost and EN placement cost.

\subsection{Heuristic algorithm (HEU)}
\label{Rand_ALG}
This section presents the heuristic scheme (\textit{HEU}). The platform choose a subset of ENs according to demand, giving priority to areas with higher demand until the available budget is fully utilized. Algorithm \ref{Alg:HEU} summarizes the  \textit{HEU} scheme as below.

\begin{algorithm}[ht!]
\caption{Heuristic EN placement algorithm (HEU)}
\label{Alg:HEU}
\begin{algorithmic}[1]
\STATE Initialization: budget $B$, demand $\lambda_i$, network delay $d_{i,j}$
\REPEAT 
  \STATE Sort $\lambda_i$ in decreasing order.
  \STATE For each area $i$, select the nearest places $ \Tilde{j}$ and calculate the EN placement cost $\mathcal{C}^{\sf HEU} = \sum_{\Tilde{j}} f_{ \Tilde{j}} \Tilde{y}_{ \Tilde{j}}$.
\UNTIL {B - $\mathcal{C}^{\sf HEU}$ $\leq 0$}\\
\STATE 
Output: EN placement decision $(\Tilde{y}_{j})$.
\end{algorithmic}
\end{algorithm}

\subsection{Two-stage SO with uniform in-sample distribution (SO)}
This section presents a two-stage stochastic optimization with uniform in-sample distribution, neglecting decisions-dependency in the demand function. In this stochastic model, it is assumed that demand $\lambda_i$ follows a uniform distribution (i.e., $p_n = \frac{1}{N}$) , with each scenario being assigned equal weight. For each scenario $n$, let $x_{i,j}^{n}$ and $u_i$ be the amount of allocated workload from AP $i$ to EN $j$, the amount of unmet demand at AP $i$ respectively. The objective of this \textit{SO} model is to optimize the expected cost over all scenarios:
\begin{subequations}
\begin{align}
    & \!\!\min_{\by,\bx,\bu}~  \sum_{j} \! f_{j} y_{j} \!\!+\! \sum_{n = 1}^{N} p_{n} \!\bigg[ \rho\!\! \sum_{i,j} \!d_{i,j} x_{i,j}^{n} \!\!+\!\! \sum_{i} s_{i} u_{i}^{n} \bigg]\!. \!\! \label{SO_obj} \\
    & \text{s.t} ~ \sum_{j} f_j y_j \leq B; ~~ \sum_{j} y_j \geq  K^{\sf min} \\
    & ~ u_i^{n} + \sum_{j} x_{i,j}^{n}  = \lambda_i^{n}, ~ \forall i,n\\
    & 0 \leq x_{i,j}^{n} \leq C_{i,j} y_j, ~ \forall i,j,n\\
    & \sum_{j} \frac{d_{i,j} x_{i,j}^{n}}{\lambda_i^{n}} \leq \Delta_i, ~ \forall i,n\\
    & \by \in \{0,1\}^{J}, ~ \bu \in \mathbb{R}_{+}^{I \times N}, ~ \bx \in \mathbb{R}_{+}^{I \times J \times N}
\end{align}
\end{subequations}
This resulting problem is a large scale mixed integer linear programming problem (MILP), which can be solved by well-known Sample average approximation (SAA) algorithm.

\subsection{Two-stage DRO with exogenous demand (DRO-DIU)}
This section presents a two-stage distributionally robust optimization with exogenous stochastic demand, neglecting decisions-dependency in the demand function. The resulting problem can be expressed as:
\begin{subequations}
\label{dro_diu}
\begin{align}
    & \!\!(\mathcal{P}_1) ~ \min_{\by} \!\sum_{j}\! f_j y_j \!+\!\! \max_{p \in \mathcal{U}} \!\min_{\bx, \bu} \mathbb{E}_{p}\!\bigg[\!\rho \! \sum_{i,j}\! d_{i,j} x_{i,j} \!+\!\! \sum_{i}\! s_i u_i\! \bigg] \!\!\label{dro_obj}\\
    &  \text{s.t.} ~~  (\ref{budget}),  (\ref{reliability}),  (\ref{var_constr1})  \nonumber \\
    &  \qquad \Omega_2 (y,\lambda) = \bigg\{ ~ 0 \leq x_{i,j} \leq C_{i,j} y_j, ~ \forall i,j \\
    & \qquad\qquad\qquad\qquad  u_i + \sum_{j} x_{i,j} = \lambda_i , ~ \forall i \\
    & \qquad\qquad\qquad\qquad \sum_{j} d_{i,j} x_{i,j} \leq \Delta_{i} \lambda_i , ~ \forall i \bigg\}.
\end{align}
\end{subequations}
The \textit{DRO-DIU} problem remains a trilevel min-max-min optimization problem, similar to the \textit{DRO-DDU} model, but with a key difference: in \textit{DRO-DIU}, the placement decision does not influence demand, leading to the setting of impact factors at zero. This issue can also be addressed using an existing algorithm by merely assigning $\Psi_{i,j}^{\mu} \!=\! \Psi_{i,j}^{\sigma} \!=\! 0, \forall i,j$.

\end{document}